\documentclass[12pt]{article}

\nonstopmode

\usepackage[left=2cm,right=2cm,top=2cm,bottom=2cm]{geometry}
\usepackage[utf8]{inputenc}
\usepackage[english]{babel}
\usepackage[T1]{fontenc}
\usepackage{amsmath}
\usepackage{amsthm}
\usepackage{amsfonts}
\usepackage{amssymb}
\usepackage{color}
\usepackage{multicol}
\usepackage{hyperref}
\hypersetup{
colorlinks=true, 
breaklinks=true, 
urlcolor= blue, 
linkcolor= blue, 
citecolor=blue, 
}
\usepackage{graphicx}
\usepackage{lscape}
\usepackage{setspace}

\newcommand{\norm}[1]{\left\lVert #1 \right\rVert}

\newtheorem{theorem}{Theorem}[section]
\newtheorem{lemma}{Lemma}[section]

\newtheorem{con}{Conjecture}[section]

\newtheorem{definition}{Definition}[section]

\title{On solutions of the Diophantine equation   \\ $F_{n_{1}}+F_{n_{2}}+F_{n_{3}}+F_{n_{4}}=2^a$ } 

\author{Pagdame TIEBEKABE $\&$ Ismaïla DIOUF\\}
\date{\today}
\begin{document}
\maketitle
\begin{abstract}
Let $(F_n)_{n\geq 0}$ be the Fibonacci sequence given by $F_0 = 0, F_1 = 1$ and $F_{n+2} = F_{n+1}+F_n$ for $n \geq 0$. In this paper, we solve all powers of two which are sums of four Fibonacci numbers with a few exceptions that we characterize. 
\end{abstract}

\section{Introduction} 
The equation $F_n + F_m=y^a$, contrary to $F_n-F_m =y^a$ has been less studied. For instance Z. Siar and R. Keskin \cite{6} have found all the solutions for $y=2$, B. Demirtürk et al \cite{7} and P. Tiebekabe et al \cite{16} have  independently determined all the solutions for y=3 and finally  F. Erduvan and R. Keskin \cite{8} for the case $y=5$ and they presumed that there are no solutions for $y>7$ prime.
Concerning the works on the the sum of Fibonacci numbers, the credit goes to Bravo and Luca \cite{9} who were the first to tackle this equation by resolving the equation $F_n+F_m=2^a$.
The result has been generalized by Pink and Ziegler in \cite{10}. In 2015, J. J Bravo and E. Bravo have determined in \cite{5} the powers of $2$ as the sums of three Fibonacci numbers and in the comments, they said that they hope it is possible to come to the end of the equation 

\begin{equation}\label{eq:1}
F_{n_{1}}+F_{n_{2}}+F_{n_{3}}+F_{n_{4}}=2^a
\end{equation}
in positive integers $n_1, n_2, n_3, n_4\quad \text{and}\quad a$ with $n_1\geqslant n_2\geqslant n_3\geqslant n_4$,
using the same method.
But since the number of cases to be considered is difficult, they leave it for other researchers.
It is in that idea that we have decided to tackle this difficult case.

Another similar problems as the one discussed in this paper were investigated for the Fibonacci and Lucas sequences. For example, repdigits which are sums of at most three Fibonacci numbers were found in \cite{12}; repdigits as sums of four Fibonacci or Lucas numbers  were found in \cite{15}; Fibonacci numbers which are sums of two repdigits were obtained in \cite{14}, while factorials which are sums of at most three Fibonacci numbers were found in \cite{13}.

After a deep analysis of works on powers of $2,3$ and $5$ as sums/differences
of Lucas or Fibonacci numbers, we think that it would be interesting to
generalize all perfect powers which are sums or differences of linear re-
currences. We should admit that recent work suggests that this will be
diffcult and may not happen any time soon. When solving $F_n \pm F_m =
y^a$ in \cite{11}, the authors were compelled to use the $abc-$conjecture to show that the equation admits finitely-many integer solutions $(n, m, y, a)$
with $\min \{y,a\}\geqslant 2$.

Let’s recall that this conjecture is still an open problem and so is not yet resolved.
Recall that the Zeckendorf representation \cite{1} of a positive integer N is the representation
$$N = F_{m_{1}} + F_{m_{2}} + \ldots + F_{m_{t}} ;\quad \text{with} \quad m_i-m_{i+1}\geqslant 2\quad \text{for}\quad  i=1,\ldots, t-1.$$
Equation (\ref{eq:1}) is a particular case of Zeckendorf representation with $N=2^a$ and $t=4$.

This paper is subdivided as follows: In Section $2$, we introduce auxyliary results used in Section $3$ to prove the main theorem of this paper stated below. 

\begin{theorem}\label{theo4}
All non-trivial solutions of the Diophantine equation (\ref{eq:1}) in positive integers $n_1, n_2, n_3, n_4$ and $a $ with $n_1\geqslant n_2\geqslant n_3\geqslant n_4$ are:

\tiny{
\begin{tabular}{|l|l|l|l|l|l|l|l|l|}\hline
$F_5+3F_2=2^3$&  $F_7+3F_2=2^4$ & $2F_4+2F_2=2^3$& $F_{13}+F_8+2F_2=2^8$\\ 
\hline
 $F_4+2F_3+F_2=2^3$&$F_6+F_5+F_3+F_2=2^4$ & $F_8+F_6+F_3+F_2= 2^5$ & $F_{16}+F_9+F_3+F_2= 2^{10}$ \\ \hline
 $F_{10}=F_5+F_4+F_2=2^6$ & $3F_5+F_2=2^4$& $F_8+2F_5+F_2=2^5$& $2F_7+F_5+F_2=2^5$\\ 
\hline
$F_9+F_8+F_6+F_2=2^6$&  $3F_8+F_2=2^6$& $F_{10}+F_5+2F_3=2^6$& $F_6+2F_4+F_3=2^4$   \\ 
\hline 
$F_{11}+F_9+F_4+F_3=2^7$& $F_{13}+F_7+F_6+F_3=2^8$& $F_{12}+F_{11}+F_8+F_3=2^8$& $F_{12}+2F_{10}+F_3=2^8$ \\ 
\hline
$F_{10}+3F_4=2^6$& $2F_{5}+2F_4=2^4$ & $F_8+F_5+2F_4=2^5$ &$2F_7+2F_4=2^5$ \\ 
\hline
$F_7+2F_6+F_4=2^5$& $F_{16}+F_8+F_7+F_4=2^{10}$&  $F_{15}+F_{14}+F_9+F_4=2^{10}$& $F_{13}+F_7+2F_5=2^8$\\ 
\hline
$F_{11}+F_8+F_7+F_5=2^7$& $2F_{10}+F_7+F_5=2^7$ &$F_{10}+2F_9+F_5=2^7$ & $F_{16}+F_8+2F_6=2^{10}$ \\ \hline
$F_{11}+3F_7=2^7$& $4F_6=2^5$ & & \\ \hline
\end{tabular}
}
\end{theorem}

\section{Auxiliary results}
In this section, we give some important known definitions, proprieties, theorem and lemmas. 
\begin{definition}

 For all algebraic numbers  $\gamma$, we define its measure by the following identity :
\begin{equation*}
{\rm M}(\gamma)=|a_d|\prod\limits_{i=1}^d \max \{1,|\gamma_{i}|\},
\end{equation*}
where $\gamma_{i}$ are the roots of $f(x)=a_d\prod\limits_{i=1}^d(x-\gamma_{i})$ is the minimal polynomial of  $\gamma$.
\end{definition} 
 
Let us define now another height, deduced from the last one, called the absolute logarithmic height. 
\begin{definition}( Absolute logarithmic height)

For a non-zero algebraic number of degree $d$ on $\mathbb{Q}$ where the minimal polynomial on  $\mathbb{Z}$ is $f(x)=a_d\prod\limits_{i=1}^d(x-\gamma_{i})$, we denote by
\begin{equation*}
h(\gamma)=\dfrac{1}{d}\left(\log|a_d|+\sum\limits_{i=1}^d \log\max\{1,|\gamma_i|\}\right) = \dfrac{1}{d} \log {\rm M}(\gamma).
\end{equation*}
the usual logarithmic absolute height of $\gamma$.
 \end{definition}
The following properties of the logarithmic height are well-known:
 \begin{itemize}
 \item[•] $h(\gamma \pm \eta)\leq h(\gamma)+h(\eta)+\log 2$;
 \item[•] $h(\gamma\eta^{\pm 1})\leq h(\gamma)+h(\eta)$;
 \item[•] $h(\gamma^k)=|k|h(\gamma)\quad k\in \mathbb{Z}$.
 \end{itemize}
The $n$th Fibonacci number can be represented as
$$
F_n=\dfrac{\alpha^n-\beta^n}{\sqrt{5}}\quad \text{for all}\quad n\geqslant 0.
$$
where $(\alpha,\beta):=((1+\sqrt{5})/2,(1-\sqrt{5})/2)$. The following inequalities 
$$
\alpha^{n-2}\leqslant F_n \leqslant \alpha^{n-1}
$$
are well-known to hold for all $n\geqslant 1$ and can be proved by induction on $n$.
The following theorem is deduced from Corollary $2.3$ of Matveev \cite{3}.
\begin{theorem} (Matveev \cite{3})\label{theo3}

 Let $n\geq 1$ an integer. Let $\mathbb{L}$ be  a field of algebraic number of degree $D$. Let $\eta_1$, \dots, $\eta_l$ non-zero elements of \ $\mathbb{L}$ and let 
$b_1$, $b_2$, \dots, $b_l$ integers, 
$$
B:=\max\{|b_1|,...,|b_l|\},
$$ 

and
\begin{equation*}
\Lambda:=\eta_1^{b_1}\cdots\eta_l^{b_l}-1=\left(\prod\limits_{i=1}^l \eta_i^{b_i}\right)-1.
\end{equation*}
Let $A_1$, \dots, $A_l$ reals numbers such that 
\begin{equation*}
A_j\geq \max\{Dh(\eta_j),|\log (\eta_j)|,0.16\}, 1\leq j\leq l.
\end{equation*}
Assume that  $\Lambda\neq 0$, So we have
\begin{equation*}
\log|\Lambda|>-3\times 30^{l+4}\times (l+1)^{5.5}\times d^2 \times A_1...A_l(1+\log D)(1+\log nB).
\end{equation*}
Further, if\  $\mathbb{L}$ is real, then 
\begin{equation*}
\log|\Lambda|>-1.4\times 30^{l+3}\times (l)^{4.5}\times d^2 \times A_1...A_l(1+\log D)(1+\log B).
\end{equation*}

\end{theorem}
The two following Lemmas are due Dujella and Peth\H o, and to Legendre respectively.

For a real number $X$, we write $\norm{X} :=\min\{\mid X-n\mid :n\in \mathbb{Z}\}$ for the distance of $X$ to the nearest integer.

\begin{lemma}\label{lem1}(Dujella and Peth\H o, \cite{2}) 

Let $M$ a positive integer, let $p/q$ the convergent of the continued fraction expansion of $\kappa$ such that $q>6M$ and let  $A$, $B$, $\mu$ real numbers such that  $A>0$ and $B>1$. Let $\varepsilon:=\norm{\mu q}-M\norm{\kappa q}$.\\ If $\varepsilon>0$ then  there is no solution of the inequality 
\begin{equation*}
0<m\kappa -n+\mu< AB^{-m}
\end{equation*}
 in integers $m$ and $n$ with 
$$
\dfrac{\log (Aq/\varepsilon)}{\log B}\leqslant m \leqslant M.
$$

\end{lemma}

\begin{lemma} \label{lem2}(Legendre) 

Let $\tau$ real number such that $x$, $y$ are integers such that 
\begin{equation*}
\left| \tau - \dfrac{x}{y}\right|<\dfrac{1}{2y^{2}},
\end{equation*}
  then $\dfrac{x}{y}=\dfrac{p_{k}}{q_{k}}$  is the convergence of $\tau$. 
  
\end{lemma}
Further, 
  \begin{equation*}
\left| \tau - \dfrac{x}{y}\right|>\dfrac{1}{(q_{k+1}+2)y^{2}}.
\end{equation*}
\section{Main result}
\begin{proof}
Assume that
\begin{equation*}
F_{n_{1}}+F_{n_{2}}+F_{n_{3}}+F_{n_{4}}=2^a
\end{equation*}
holds.\\
Let us first find relation between $n_1$ and $a$. 

Combining equation(\ref{eq:1}) with the well-known inequality $F_n\leqslant \alpha^{n-1}$ for all $n\geqslant 1$, one gets that 
\begin{align*}
F_{n_{1}}+F_{n_{2}}+F_{n_{3}}+F_{n_{4}}=2^a\leqslant& \alpha^{n_1-1}+\alpha^{n_2-1}+\alpha^{n_3-1}+\alpha^{n_4-1}\\
<&2^{n_1-1}+2^{n_2-1}+2^{n_3-1}+2^{n_4-1}\quad \because \alpha<2\\
<& 2^{n_1-1}\left(1+2^{n_2-n_1}+2^{n_3-n_1}+2^{n_4-n_1}\right)\\
\leqslant & 2^{n_1-1}\left( 1+1+2^{-1}+2^{-2} \right)= 2^{n_1-1}\left(2+2^{-1}+2^{-2}\right)\\
<& 2^{n_1+1}.
\end{align*}
Hence
$$
2^a<2^{n_1+1}\Longrightarrow a< n_1+1\Longrightarrow a\leqslant n_1.
$$
This inequality will help  us to calculate some parameters.

Rewriting equation (\ref{eq:1}), we get
$$
\dfrac{\alpha^{n_1}}{\sqrt{5}}-2^a= \dfrac{\beta^{n_1}}{\sqrt{5}}-(F_{n_{2}}+F_{n_{3}}+F_{n_{4}}).
$$
Taking absolute values on the above equation, we obtain
$$
\left|\dfrac{\alpha^{n_1}}{\sqrt{5}}-2^a \right|\leqslant \left| \dfrac{\beta^{n_1}}{\sqrt{5}}\right|+ (F_{n_{2}}+F_{n_{3}}+F_{n_{4}})< \dfrac{|\beta|^{n_1}}{\sqrt{5}}+ (\alpha^{n_{2}}+\alpha^{n_{3}}+\alpha^{n_{4}}),
$$
and
$$
\left|\dfrac{\alpha^{n_1}}{\sqrt{5}}-2^a \right|< \dfrac{1}{2}+ \left(\alpha^{n_{2}}+\alpha^{n_{3}}+\alpha^{n_{4}}\right)\quad \text{where we used}\quad F_n\leqslant \alpha^{n-1}.
$$

Dividing both side of the above equation by $\alpha^{n_{1}}/\sqrt{5}$, we get 

$$
\left|1-2^a\cdot \alpha^{-n_1}\cdot \sqrt{5}\right|< \dfrac{\sqrt{5}}{2\alpha^{n_1}}+ \left(\alpha^{n_{2}-n_1}+\alpha^{n_{3}-n_1}+\alpha^{n_{4}-n_1}\right)\sqrt{5}< \dfrac{\sqrt{5}}{2\alpha^{n_1}}+\dfrac{\sqrt{5}}{\alpha^{n_1-n_2}}+\dfrac{\sqrt{5}}{\alpha^{n_1-n_3}}+\dfrac{\sqrt{5}}{\alpha^{n_1-n_4}}.
$$

Taking account the assumption $n_4\leqslant n_3\leqslant n_3\leqslant n_2\leqslant n_1$, we get

\begin{equation}\label{eq:2}
|\Lambda_1|=\left|1-2^a\cdot \alpha^{-n_1}\cdot \sqrt{5}\right|< \dfrac{9}{\alpha^{n_1-n_2}}
\end{equation}

Let apply Matveev's theorem, with the following parameters $t:=3$ and 
$$
\gamma_1:=2,\quad \gamma_2:=\alpha,\quad \gamma_3:=\sqrt{5},\quad b_1:=a,\quad b_2:=-n,\quad \text{and}\quad b_3:=1.
$$

Since $\gamma_1, \gamma_2, \gamma_3\in \mathbb{K}:=\mathbb{Q}(\sqrt{5}),$ we can take $D:=2$. Before to apply Matveev's theorem, we have to check the last condition: the left-hand side of (\ref{eq:2}) is not zero. Indeed, if it were zero, we would then get that $2^a\sqrt{5}=\alpha^n$. Squaring the previous relation, we get $\alpha^{2n}=5\cdot 2^{2a}=5\cdot 4^a$. This implies that $\alpha^{2n}\in \mathbb{Z}$. Which is impossible. Then $\Lambda_1\neq 0$. The logarithmic height of $\gamma_1, \gamma_2$ and $\gamma_3$ are:

$h(\gamma_1)=\log 2= 0.6931\ldots$, so we can choose $A_1:=1.4$.

$h(\gamma_2)=\dfrac{1}{2}\log \alpha=0.2406\ldots$, so we can choose $A_2:= 0.5$.

$h(\gamma_3)=\log \sqrt{5}=0.8047\ldots$, it follows that we can choose $A_3:=1.7$.

Since $a< n_1+1$, $B:=\max \{|b_1|, |b_2|, |b_3|\}=n_1.$ Matveev's result informs us that
\begin{equation}\label{eq:3}
\left|1-2^a\cdot \alpha^{n_1}\cdot \sqrt{5}\right|>\exp \left(-c_1\cdot (1+\log n)\cdot 1.4\cdot 0.5\cdot 1.7\right),
\end{equation}
where $c_1:=1.4\cdot 30^6\cdot 3^{4.5}\cdot 2^2\cdot (1+\log 2)<9.7\times 10^{11}.$

Taking $\log$ in inequality (\ref{eq:2}), we get
$$
\log |\Lambda_1|< \log 9-(n_1-n_2)\log \alpha.
$$

Taking $\log$ in inequality (\ref{eq:3}), we get
$$
\log |\Lambda_1|> 2.31\times 10^{12}\log n_1.
$$

Comparing the previous two inequality, we get
$$
(n_1-n_2)\log \alpha-\log 9< 2.31\times 10^{12}\log n_1,
$$
where we used $1+\log n_1< 2\log n_1$ which holds for all $n_1\geqslant 3$. Then we have
\begin{equation}\label{eq:4}
(n_1-n_2)\log \alpha<2.32\times 10^{12}\log n_1.
\end{equation}
Let us now consider a second linear form in logarithms. Rewriting equation (\ref{eq:1}) as follows
$$
\dfrac{\alpha^{n_1}}{\sqrt{5}}+\dfrac{\alpha^{n_2}}{\sqrt{5}}-2^a=\dfrac{\beta^{n_1}}{\sqrt{5}}+\dfrac{\beta^{n_2}}{\sqrt{5}}-\left(F_{n_3}+F_{n_4}\right).
$$
Taking absolute values on the above equation and the fact that $\beta=(1-\sqrt{5})/2$, we get
\begin{align*}
\left| \dfrac{\alpha^{n_1}}{\sqrt{5}}\left(1+\alpha^{n_2-n_1}\right) -2^a\right|\leqslant &\dfrac{|\beta|^{n_1}+|\beta|^{n_2}}{\sqrt{5}}+F_{n_3}+F_{n_4}\\
<& \dfrac{1}{3}+ \alpha^{n_3}+\alpha^{n_4}\quad \text{for all} \quad n_1\geqslant 5\quad \text{and}\quad n_2\geqslant 5.
\end{align*}
Dividing both sides of the above inequality by $\dfrac{\alpha^{n_1}}{\sqrt{5}}\left(1+\alpha^{n_2-n_1}\right)$, we obtain 
\begin{equation}\label{eq:5}
|\Lambda_2|=\left|1-2^a\cdot\alpha^{n_1}\cdot\sqrt{5}\left(1+\alpha^{n_2-n_1}\right)^{-1}\right|<\dfrac{6}{\alpha^{n_2-n_1}}.
\end{equation}
Let apply  Matveev's theorem second time with the following data 
$$
t:=3,\quad \gamma_1:=2,\quad \gamma_2:=\alpha,\quad \gamma_3:=\sqrt{5}\left(1+\alpha^{n_2-n_1}\right)^{-1},\quad b_1:=a,\quad b_2:=-n_1,\quad \text{and}\quad b_3:=1.
$$
Since $\gamma_1, \gamma_2, \gamma_3\in \mathbb{K}:=\mathbb{Q}(\sqrt{5}),$ we can take $D:=2$. The left hand side of (\ref{eq:5}) is not zero, otherwise, we would get the relation
\begin{equation}\label{eq:6}
2^a\sqrt{5}=\alpha^{n_1}+\alpha^{n_2}.
\end{equation}
Conjugating (\ref{eq:6}) in the field $\mathbb{Q}(\sqrt{5})$, we get
\begin{equation}\label{eq:7}
-2^a\sqrt{5}=\beta^{n_1}+\beta^{n_2}.
\end{equation}

Combining (\ref{eq:6}) and (\ref{eq:7}), we get 
$$
\alpha^{n_1}<\alpha^{n_1}+\alpha^{n_2}=|\beta^{n_1}+\beta^{n_2}|\leqslant |\beta|^{n_1}+|\beta|^{n_2}< 1
$$

which is impossible for $n_1>350$. Hence $\Lambda_2\neq 0$.
We know that, $h(\gamma_1)=\log 2$ and $h(\gamma_2)=\dfrac{1}{2}\log \alpha$.
Let us now estimate $h(\gamma_3)$ by first observing that
$$
\gamma_3=\dfrac{\sqrt{5}}{1+\alpha^{n_2-n_1}}<\sqrt{5}\quad \text{and}\quad \gamma_3^{-1}=\dfrac{1+\alpha^{n_2-n_1}}{\sqrt{5}}<\dfrac{2}{\sqrt{5}},
$$
so that $|\log\gamma_3|<1.$
Using proprieties of logarithmic height stated in Section $2$, we have
$$
h(\gamma_3)\leqslant \log \sqrt{5}+ |n_2-n_1|\left(\dfrac{\log \alpha}{2}\right)+\log 2= \log (2\sqrt{5})+ (n_1-n_2)\left(\dfrac{\log \alpha}{2}\right).
$$ 
Hence, we can take $A_3:=3+(n_1-n_2)\log \alpha>\max \{2h(\gamma_3), |\log \gamma_3|,0.16\}$. 

Matveev's theorem implies that 
$$
\exp \left(-c_2(1+\log n_1)\cdot 1.4\cdot 0.5\cdot(3+(n_1-n_2)\log \alpha)\right)
$$
where $c_2:=1.4\cdot 30^6\cdot 3^{4.5}\cdot 2^2\cdot (1+\log 2)<9.7\times 10^{11}.$

Since $(1+\log n_1)< 2\log n_1$ hold for $n_1\geqslant 3,$ from (\ref{eq:5}), we have
\begin{equation}\label{eq:8}
(n_1-n_3)\log \alpha-\log 6< 1.4\times 10^{12}\log n_1(3+(n_1-n_2)\log \alpha).
\end{equation} 

Putting relation (\ref{eq:4}) in the right-hand side of (\ref{eq:8}), we get

\begin{equation}\label{eq:9}
(n_1-n_3)\log \alpha< 3.29\times 10^{24}\log^2 n_1.
\end{equation}
Let us consider a third linear form in logarithms. To this end, we again rewrite (\ref{eq:1}) as follows
$$
\dfrac{\alpha^{n_1}+\alpha^{n_2}+\alpha^{n_3}}{\sqrt{5}}-2^a=\dfrac{\beta^{n_1}+\beta^{n_2}+\beta^{n_3}}{\sqrt{5}}-F_{n_4}.
$$
Taking absolute values on both sides, we obtain
\begin{align*}
\left|\dfrac{\alpha^{n_1}}{\sqrt{5}}\left(1+\alpha^{n_2-n_1}+\alpha^{n_3-n_1}\right)-2^a\right|\leqslant & \dfrac{|\beta|^{n_1}+|\beta|^{n_2}+|\beta|^{n_3}}{\sqrt{5}}+F_{n_4}\\
< & \dfrac{3}{4}+\alpha^{n_4}\quad \text{for all}\quad n_1>350,\quad \text{and}\quad n_2, n_3, n_4\geqslant 1.
\end{align*}
Thus we have
\begin{equation}\label{eq:10}
|\Lambda_3|=\left|1-2^a\cdot \alpha^{-n_1}\cdot\sqrt{5}\left(1+\alpha^{n_2-n_1}+\alpha^{n_3-n_1}\right)^{-1}\right|< \dfrac{3}{\alpha^{n_1-n_4}}.
\end{equation}

In a third application of Matveev's theorem, we can take parameters
$$
t:=3,\quad \gamma_1:=2,\quad \gamma_2:=\alpha,\quad \gamma_3:=\sqrt{5}\left(1+\alpha^{n_2-n_1}+\alpha^{n_3-n_1}\right)^{-1},\quad b_1:=a,\quad b_2:=-n,\quad \text{and},\quad b_3:=1.
$$
Since $\gamma_1, \gamma_2, \gamma_3\in \mathbb{K}:=\mathbb{Q}(\sqrt{5}),$ we can take $D:=2$. The left hand side of (\ref{eq:10}) is not zero. The proof is done by contradiction. Suppose the contrary. Then
$$
2^a\sqrt{5}=\alpha^{n_1}+\alpha^{n_2}+\alpha^{n_3}.
$$
Taking the conjugate in the field $\mathbb{Q}(\sqrt{5})$, we get
$$
-2^a\sqrt{5}=\beta^{n_1}+\beta^{n_2}+\beta^{n_3},
$$
which leads to
$$
\alpha^{n_1}< \alpha^{n_1}+\alpha^{n_2}+\alpha^{n_3}=|\beta^{n_1}+\beta^{n_2}+\beta^{n_3}|\leqslant |\beta|^{n_1}+|\beta|^{n_2}+|\beta|^{n_3}< 1
$$
and leads to a contradiction since $n_1>350.$ Hence $\Lambda_3\neq 0$.

As before we did, we can take $A_1:=1.4, A_2:=0.5$ and $B:=n_1.$ We can also see that 
$$
\gamma_3=\dfrac{\sqrt{5}}{1+\alpha^{n_2-n_1}+\alpha^{n_3-n_1}}< \sqrt{5}\quad \text{and}\quad \gamma_3^{-1}=\dfrac{1+\alpha^{n_2-n_1}+\alpha^{n_3-n_1}}{\sqrt{5}}<\dfrac{3}{\sqrt{5}},
$$
so $|\log \gamma_3|<1$. Applying proprieties on logarithmic height, we estimate $h(\gamma_3)$.
Hence 
\begin{align*}
h(\gamma_3)\leqslant & \log \sqrt{5}+|n_2-n_1|\left(\dfrac{\log\alpha}{2}\right)+|n_3-n_1|\left(\dfrac{\log\alpha}{2}\right)+\log 3\\
=& \log(3\sqrt{5})+(n_1-n_2)\left(\dfrac{\log\alpha}{2}\right)+(n_1-n_3)\left(\dfrac{\log\alpha}{2}\right);
\end{align*}
so we can take
$$
A_3:=4+(n_1-n_2)\log \alpha+(n_1-n_3)\log \alpha> \max \{2h(\gamma_3), |\log \gamma_3|, 0.16\}.
$$
A lower bound on the left-hand side of (\ref{eq:10}) is
$$
\exp (-c_3\cdot(1+\log n_1)\cdot 1.4\cdot 0.5 \cdot(4+(n_1-n_2)\log \alpha+ (n_1-n_3)\log \alpha))
$$

where $c_3=1.4\cdot 30^6\cdot 3^{4.5}\cdot 2^2\cdot (1+\log 2)<9.7\times 10^{11}.$

From inequality (\ref{eq:10}), we have 
\begin{equation}\label{eq:11}
(n_1-n_4)\log \alpha< 1.4\times 10^{12}\log n_1\cdot(4+(n_1-n_2)\log \alpha+(n_1-n_3)\log \alpha).
\end{equation}
Combining equation (\ref{eq:4}) and (\ref{eq:9}) in the right-most terms of equation (\ref{eq:11}) and performing the respective calculations, we get
\begin{equation}\label{eq:12}
(n_1-n_4)\log \alpha< 9.3\times \log^3n_1.
\end{equation} 

Let us now consider a forth and last linear form in logarithms. Rerwriting (\ref{eq:1}) once again by separating large terms and small terms, we get

$$
\dfrac{\alpha^{n_1}+\alpha^{n_2}+\alpha^{n_3}+\alpha^{n_4}}{\sqrt{5}}-2^a=\dfrac{\beta^{n_1}+\beta^{n_2}+\beta^{n_3}+\beta^{n_4}}{\sqrt{5}}.
$$
Taking absolute values on both sides, we get
\begin{align*}
\left|\dfrac{\alpha^{n_1}}{\sqrt{5}}\left(1+\alpha^{n_2-n_1}+\alpha^{n_3-n_1}+\alpha^{n_4-n_1}\right)-2^a\right|\leqslant & \dfrac{|\beta|^{n_1}+|\beta|^{n_2}+|\beta|^{n_3}+|\beta|^{n_4}}{\sqrt{5}}\\
< & \dfrac{4}{5}\quad \text{for all}\quad n_1>350,\quad \text{and}\quad n_2, n_3, n_4\geqslant 1.
\end{align*}

Dividing both sides of the above relation by the fist term of the RHS of the previous equation, we get 
\begin{equation}\label{eq:13}
|\Lambda_4|=\left|1-2^a\cdot\alpha^{-n_1}\cdot \sqrt{5}\left(1+\alpha^{n_2-n_1}+\alpha^{n_3-n_1}+\alpha^{n_4-n_1}\right)^{-1}\right|<\dfrac{2}{\alpha^{n_1}}.
\end{equation}

In the last application of Matveev's theorem, we have the following parameters
$$
\gamma_1:=2,\quad \gamma_2:=\alpha,\quad \gamma_3:=\sqrt{5}\left(1+\alpha^{n_2-n_1}+\alpha^{n_3-n_1}+\alpha^{n_4-n_1}\right)^{-1},
$$
and we can also take $b_1:=a,\quad b_2:=-n$ and $b_3:=1.$ Since $\gamma_1, \gamma_2, \gamma_3\in \mathbb{K}:=\mathbb{Q}(\sqrt{5}),$ we can take $D:=2$. The left hand side of (\ref{eq:13}) is not zero. The proof is done by contradiction. Suppose the contrary. Then
$$
2^a\sqrt{5}=\alpha^{n_1}+\alpha^{n_2}+\alpha^{n_3}+\alpha^{n_4}.
$$
Conjugating the above relation in the field $\mathbb{Q}(\sqrt{5})$, we get
$$
-2^a\sqrt{5}=\beta^{n_1}+\beta^{n_2}+\beta^{n_3}+\beta^{n_4}.
$$
Combining the above two equations, we get
$$
\alpha^{n_1}< \alpha^{n_1}+\alpha^{n_2}+\alpha^{n_3}+\alpha^{n_4}=|\beta^{n_1}+\beta^{n_2}+\beta^{n_3}+\beta^{n_4}|\leqslant |\beta|^{n_1}+|\beta|^{n_2}+|\beta|^{n_3}+|\beta|^{n_4}< 1,
$$
and leads to contradiction since $n_1>350.$

As before, here, we can take $A_1:=1.4, A_2:=0.5$ and $B:=n_1$. Let us estimate $h(\gamma_3)$.
We can see that, 
$$
\gamma_3=\dfrac{\sqrt{5}}{1+\alpha^{n_2-n_1}+\alpha^{n_3-n_1}+\alpha^{n_4-n_1}}<\sqrt{5}\quad \text{and}\quad \gamma_3^{-1}=\dfrac{1+\alpha^{n_2-n_1}+\alpha^{n_3-n_1}+\alpha^{n_4-n_1}}{\sqrt{5}}<\dfrac{4}{\sqrt{5}}.
$$
Hence $|\log \gamma_3|<1.$ Then 
\begin{align*}
h(\gamma_3)\leqslant &\log (4\sqrt{5})+|n_2-n_1|\left(\dfrac{\log\alpha}{2}\right)+|n_3-n_1|\left(\dfrac{\log\alpha}{2}\right)+|n_4-n_1|\left(\dfrac{\log\alpha}{2}\right)\\
=& \log(4\sqrt{5})+(n_1-n_2)\left(\dfrac{\log\alpha}{2}\right)+(n_1-n_3)\left(\dfrac{\log\alpha}{2}\right)+(n_1-n_4)\left(\dfrac{\log\alpha}{2}\right);
\end{align*}
so we can take
$$
A_3:= 5+(n_1-n_2)\log \alpha+(n_1-n_3)\log \alpha+ (n_1-n_4)\log \alpha.
$$
Then a lower bound on the left-hand side  of (\ref{eq:13}) is 
$$
\exp (-c_4\cdot(1+\log n_1)\cdot 1.4\cdot 0.5\cdot(5+(n_1-n_2)\log \alpha+(n_1-n_3)\log \alpha+ (n_1-n_4)\log \alpha)),
$$
where $c_4=1.4\cdot 30^6\cdot 3^{4.5}\cdot 2^2\cdot (1+\log 2)<9.7\times 10^{11}.$

So, inequality (\ref{eq:13}) yields 
\begin{equation}\label{eq:14}
n_1\log \alpha< 1.4\times 10^{12}\log n_1\cdot (5+(n_1-n_2)\log \alpha+(n_1-n_3)\log \alpha+ (n_1-n_4)\log \alpha).
\end{equation}
Using now (\ref{eq:4}), (\ref{eq:9}) and (\ref{eq:12}) in the right-most terms of the above inequality (\ref{eq:14}) and performing the respective calculation, we find that 
$$
n_1\log \alpha< 40.32\times 10^{48}\log^4n_1.
$$
With the help of \textit{Mathematica}, we get from the previous inequality 

$$
n<2.8\times 10^{58}.
$$

We record what we have proved.
\begin{lemma}
If $(n_1, n_2, n_3, n_4, a)$ is a positive solution of (\ref{eq:1}) with $n_1\geqslant n_2\geqslant n_3\geqslant n_4$, then 
$$
a\leqslant n_1< 2.8\times 10^{58}.
$$
\end{lemma}

\section{Reduction the bound on $n$}
The goal of this section is reduce the upper bound on $n$ to a size that can be handled. To do this, we shall use Lemma \ref{lem1} four times. Let us consider 
\begin{equation}\label{eq:15}
z_1:=a\log 2-n_1\log \alpha+ \log \sqrt{5}.
\end{equation}
From equation (\ref{eq:15}), (\ref{eq:2}) can be written as
\begin{equation}\label{eq:16}
\left|1-e^{z_1}\right|< \dfrac{9}{\alpha^{n_1-n_2}}.
\end{equation}
Associating (\ref{eq:1}) and Binet's formula for the Fibonacci sequence, we have
$$
\dfrac{\alpha^{n_1}}{\sqrt{5}}=F_{n_1}+\dfrac{\beta^{n_1}}{\sqrt{5}}<F_{n_{1}}+F_{n_{2}}+F_{n_{3}}+F_{n_{4}} =2^a,
$$
hence 
$$
\dfrac{\alpha^{n_1}}{\sqrt{5}}< 2^a,
$$

which leads to $z_1>0.$ This result together with (\ref{eq:16}), give
$$
0< z_1< e^{z_1}-1< \dfrac{9}{\alpha^{n_1-n_2}}.
$$
Replacing (\ref{eq:15}) in the inequality and dividing both sides of the resulting inequality by $\log \alpha$, we get
\begin{equation}\label{eq:17}
0< a \left(\dfrac{\log 2}{\log \alpha}\right)-n+ \left(\dfrac{\log \sqrt{5}}{\log \alpha}\right)<\dfrac{9}{\log \alpha}\cdot \alpha^{n_1-n_2}< 19\cdot \alpha^{n_1-n_2}.
\end{equation} 

We put 
$$
\tau:=\dfrac{\log 2}{\log \alpha}, \quad \mu:=\dfrac{\log \sqrt{5}}{\log \alpha},\quad A:=19,\quad \text{and}\quad B:=\alpha.
$$
$\tau$ is an irrational number. We also put $M:=2.8\times 10^{58}$, which is an upper bound on $a$ by Lemma \ref{lem1} applied to inequality, that 
$$
n_1-n_2< \dfrac{\log(Aq/\varepsilon)}{\log B},
$$
where $q> 6M$ is a denominator of a convergent of the continued fraction of $\tau$ such that $\varepsilon:=\norm{\mu q}-M\norm{\tau q}>0.$ A computation with \textit{SageMath} revealed that if  $(n_1, n_2, n_3, n_4, a)$ is a possible solution of the equation \ref{eq:1}, then 

$$
n_1-n_2\in [0, 314].
$$

Let us now consider a second function, derived from (\ref{eq:5}) in order to find an improved upper bound on $n_1-n_2$.

Put 
$$
z_2:=a\log 2-n_1\log \alpha+ \log \Upsilon_1(n_1-n_2)
$$
where $\Upsilon$ is the function given by the formula $\Upsilon(t):=\sqrt{5}\left(1+\alpha^{-t}\right)^{-1}.$ From (\ref{eq:5}), we have 
\begin{equation}\label{eq:18}
\left| 1-e^{z_2}\right|< \dfrac{6}{\alpha^{n_1-n_3}}.
\end{equation}
Using (\ref{eq:1}) and the Binet's formula for the Fibonacci sequence, we have
$$
\dfrac{\alpha^{n_1}}{\sqrt{5}}+\dfrac{\alpha^{n_2}}{\sqrt{5}}= F_{n_1}+F_{n_2}+\dfrac{\beta^{n_1}}{\sqrt{5}}+\dfrac{\beta^{n_2}}{\sqrt{5}}<F_{n_1}+F_{n_2}+1\leqslant F_{n_{1}}+F_{n_{2}}+F_{n_{3}}+F_{n_{4}}=2^a.
$$
Therefore $1< 2^a\sqrt{5}\alpha^{-n_1}\left(1+\alpha^{n_2-n_1}\right)^{-1}$ and so $z_2>0$. This with (\ref{eq:18}) give
$$
0<z_2\leqslant e^{z_2}-1< \dfrac{6}{\alpha^{n_1-n_3}}.
$$
Putting expression of $z_2$ in the above inequality and arguing as in (\ref{eq:17}), we obtain
\begin{equation}\label{eq:19}
0<a\left(\dfrac{\log 2}{\log \alpha}\right)-n_1+\dfrac{\log\Upsilon_1(n_1-n_2)}{\log \alpha}< 13\cdot \alpha^{-(n_1-n_3)}.
\end{equation}
As before, we take again $M:=2.8\times 10^{58}$ which is the upper bound on $a$, and, as explained before, we apply Lemma \ref{lem1} to inequamity (\ref{eq:19}) for all choices $n_1-n_2\in [0,314]$ except when $n_1-n_2=2,6$. With the help of \textit{SageMath}, we find that if $(n_1,n_2,n_3,n_4,a)$ is a possible solution of the equation (\ref{eq:1}) with $n_1-n_2\neq 2$ and $n_1-n_2\neq 6$, then $n_1-n_3\in [0, 314]$.

Study of the cases $n_1-n_2\in \{2,6\}$. For these cases, when we apply Lemma \ref{lem1} to the expression (\ref{eq:19}), the corresponding parameter $\mu$ appearing in Lemma \ref{lem1} is 
$$
\dfrac{\log \Upsilon_1(t)}{\log\alpha}=
\left\{
    \begin{array}{lllll}
        1 & \text{if} & t=2;\\
        3-\dfrac{\log2}{\log\alpha} & \text{if} & t=6.
    \end{array}
\right.
$$
In both case, the parameters $\tau$ and $\mu$ are linearly dependent, which yields that the corresponding value of $\varepsilon$ from Lemma \ref{lem1} is always negative and therefore the reduction method is not useful  for reducing the bound on $n$ in these instances. For this, we need to treat these case differently. 

However, we can see that if $t=2$ and $6$, then resulting inequality from (\ref{eq:19}) has the shape $0<|x\tau-y|<13\cdot \alpha^{-(n_1-n_3)}$ with $\tau$ being irrational number and $x,y\in \mathbb{Z}$. Then, using the known proprieties of the convergents of the continued fractions to obtain a nontrivial lower bound for $|x\tau-y|$. Let see how to do.
 
For $n_1-n_2=2$, from (\ref{eq:19}), we get that
\begin{equation}\label{eq:20}
0<a\tau -(n_1-1)< 13\cdot \alpha^{-(n_1-n_3)}, \quad \text{where}\quad \tau=\dfrac{\log 2}{\log \alpha}.
\end{equation}
Let $[a_1, a_2, a_3, a_4,\ldots]=[1,2,3,1,\ldots]$ be the continued fraction of $\tau$, and let denote $p_k/q_k$ it $k$th convergent. By Lemma \ref{lem2}, we know that $a<2.8\times 10^{58}$. An inspection in \textit{SageMath} reveals that 
$$
1207471144047491451512110092657730332808809199105354185685=q_{113}<2.8\times 10^{58}<
$$
$$
q_{114}=28351096929195187169517686575841899309129196859170938821667.
$$
Furthermore, $a_M:=\max \{a_i: i=0,1,\ldots,114\}=134.$ So, from the proprieties of continued fractions, we obtain that 
\begin{equation}\label{eq:21}
|a\tau-(n_1-1)|> \dfrac{1}{(a_M+2)a}.
\end{equation} 
Comparing (\ref{eq:20}) and (\ref{eq:21}), we get
$$
\alpha^{n_1-n_3}<13\cdot (134+2)a.
$$
Taking $\log$ on both sides of above equation and the divide the obtained result by $\log\alpha$, we get 
$$
n_1-n_3<296.
$$

In order to avoid repetition, we freely omits the details for the case $n_1-n_2=6.$ Here, we get $n_1-n_3<314$.

This completes the analysis of the two special cases $n_1-n_2=2$ and $n_1-n_2=6.$ Consequently $n_1-n_3\leqslant 314$ always holds.

Now let us use (\ref{eq:10}) in order to find improved upper bound on $n_1-n_4.$ Put 
$$
z_3:= a\log 2-n_1\log\alpha+\log\Upsilon_2(n_1-n_2, n_1-n_3),
$$ 
where $\Upsilon_2$ is the function given by the formula $\Upsilon_2(t,s):=\sqrt{5}\left(1+\alpha^{-t}+\alpha^{-s}\right)^{-1}.$ From (\ref{eq:10}), we have
\begin{equation}\label{eq:22}
|1-e^{z_3}|<\dfrac{3}{\alpha^{n_1-n_4}}.
\end{equation}

Note that, $z_3\neq 0$; thus, two cases arises: $z_3>0$ and $z_3< 0$.

If $z_3>0,$ then
$$
0<z_3\leqslant e^{z_3}-1<\dfrac{3}{\alpha^{n_1-n_4}}.
$$
Suppose now $z_3< 0.$ It is easy to check that $3/\alpha^{n_1-n_4}<1/2$ for all $n_1>350$ and $n_4\geqslant 2$. From (\ref{eq:22}), we have that 
$$
|1-e^{z_3}|<1/2 \quad \text{and therefore}\quad e^{|z_3|}<2.
$$
Since $z_3<0$, we have:
$$
0<|z_3|\leqslant e^{|z_3|}-1=e^{|z_3|}\left|e^{|z_3|}-1\right|< \dfrac{6}{\alpha^{n_1-n_4}}
$$
which give 
$$
0<|z_3|<\dfrac{6}{\alpha^{n_1-n_4}}
$$
holds for $z_3<0$,  $z_3>0$ and for all for all $n_1>350$, and $n_4\geqslant 2$.
Replacing the expression of $z_3$ in the above inequality and arguing again as before, we conclude that
\begin{equation}\label{eq:23}
0< \left|a\left(\dfrac{\log 2}{\log \alpha}\right)-n_1+ \dfrac{\log \Upsilon_2(n_1-n_2, n_1-n_3)}{\log \alpha}\right|<13\cdot \alpha^{-(n_1-n_4)}.
\end{equation}

Here, we also take, $M:=2.8\times 10^{58}$ and we apply Lemma \ref{lem1} in inequality (\ref{eq:23}) for all choices $n_1-n_2\in \{0,314\}$ and $n_1-n_3\in \{0, 314\}$ except when 
$$
(n_1-n_2,n_1-n_3)\in \{(0,3), (1,1), (1,5), (3,0), (3,4), (4,3), (5,1), (7,8), (8,7)\}.
$$
Indeed, with the help of \textit{SageMath} we find that if $(n_1, n_2, n_3, n_4, a)$ is a possible solution of the equation (\ref{eq:1}) excluding theses cases presented before. Then $n_1-n_4\leqslant 314.$

\underline{SPECIAL CASES}. We deal with the cases when 
$$
(n_1-n_2,n_1-n_3)\in \{(1,1), (3,0), (4,3), (5,1), (8,7)\}.
$$
It is easy to check that

$$
\dfrac{\log \Upsilon_2(t,s)}{\log\alpha}=
\left\{
    \begin{array}{lllll}
        0, & \text{if} & (t,s)=(1,1);\\
        0, & \text{if} & (t,s)=(3,0);\\
        1, & \text{if} & (t,s)=(4,3);\\
        2-\dfrac{\log 2}{\log \alpha}, & \text{if} & (t,s)=(5,1);\\
        3-\dfrac{\log 2}{\log \alpha}, & \text{if} & (t,s)=(8,7).
    \end{array}
\right.
$$
As we explained before, when we apply Lemma \ref{lem1} to the expression (\ref{eq:23}), the parameters $\tau$ and $\mu$ are linearly dependent, so the corresponding value of $\varepsilon$ from Lemma \ref{lem1} is always negative in all cases. For this reason, we shall treat these cases differently.

Here, we have to solve the equations
$$
F_{n_2+1}+2F_{n_2}+F_{n_4}=2^a,\quad  2F_{n_2+3}+F_{n_2}+F_{n_4}=2^a,\quad F_{n_2+4}+F_{n_2}+F_{n_2+1}+F_{n_4}=2^a,
$$
\begin{equation}\label{eq:24}
 F_{n_2+5}+F_{n_2}+F_{n_2+4}+F_{n_4}=2^a,\quad \text{and}\quad F_{n_2+8}+F_{n_2}+F_{n_2+1}+F_{n_4}=2^a
\end{equation}
in positive integers $n_2, n_4$ and $a$. To do so, we recall the following well-known relation between
the Fibonacci and the Lucas numbers:
\begin{equation}\label{eq:25}
L_k=F_{k-1}+F_{k+1} \quad \text{for all} \quad k\geqslant 1.
\end{equation}
From (\ref{eq:25}) and (\ref{eq:24}), we have the following identities 

$$
F_{n_2+1}+2F_{n_2}+F_{n_4}=F_{n_2+2}+F_{n_2}+F_{n_4}=F_{k+2}+F_{k}+F_{m}, 
$$
$$
2F_{n_2+3}+F_{n_2}+F_{n_4}=F_{n_2+2}+F_{n_2+4}+F_{n_4}=F_{k+2}+F_{k+4}+F_{m},
$$
\begin{equation}\label{eq:26}
 F_{n_2+4}+F_{n_2}+F_{n_2+1}+F_{n_4}=F_{n_2+2}+F_{n_2+4}+F_{n_4}=F_{k+2}+F_{k+4}+F_{m}, 
\end{equation}
$$
F_{n_2+5}+F_{n_2}+F_{n_2+4}+F_{n_4}=2F_{n_2+2}+2F_{n_2+4}+F_{n_4}=2F_{k+2}+2F_{k+4}+F_{m},
$$
$$
\text{and}\quad F_{n_2+8}+F_{n_2}+F_{n_2+1}+F_{n_4}= 2F_{n_2+6}+2F_{n_2+4}+F_{n_4}=2F_{k+6}+2F_{k+4}+F_{m},
$$
hold for all $k, m\geqslant 0$.

Equation (\ref{eq:24}) are transformed into the equations 
\begin{equation}\label{eq:27}
L_{k+1}+F_m=2^a,\quad L_{k+3}+F_{m}=2^a,\quad 2L_{k+3}+F_m=2^a,\quad 2L_{k+5}+F_m=2^a,
\end{equation}
to be resolved in positive integers $k,m$ and $a$.

A quick search in SageMath and analytical resolution leads to : 

$$
(k,m,a)\in \{(4,5,4), (4,8,5)\}\quad \text{for}\quad L_{k+1}+F_m=2^a,
$$

$$
(k,m,a)\in \{(2, 5, 4), (2,8,5), (4,4,5)\} \quad \text{for}\quad L_{k+3}+F_{m}=2^a,
$$

$$
(k,m,a)=(5,9,7) \quad \text{for}\quad 2L_{k+3}+F_m=2^a,
$$

$$
(k,m,a)=(3,9,7) \quad \text{for}\quad 2L_{k+5}+F_m=2^a.
$$
 A completes  resolution and analysis gives solutions that are already listed in Theorem \ref{theo4}. This completes the analysis of the special cases.
 
Finally let us use (\ref{eq:13}) in order to find improved upper bound on $n_1.$ Put
$$
z_4:=a\log 2-n_1\log \alpha+ \log\Upsilon_3(n_1-n_2, n_1-n_3, n_1-n_4),
$$ 
where $\Upsilon_3$ is the function given by the formula 
$$
\Upsilon_3(t,u,v):= \sqrt{5}\left( 1+ \alpha^{-t}+\alpha^{-u}+\alpha^{-v}\right)^{-1}
$$
with $t=n_1-n_2, u=n_1-n_3$ and $v=n_1-n_4.$ From (\ref{eq:13}), we get 
\begin{equation}\label{eq:28}
|1-e^{z_3}|<\dfrac{2}{\alpha^{n_1}}.
\end{equation}

Since $z_3\neq 0$, as before, two cases arise: $z_4< 0$ and $z_4> 0$.

If $z_4>0$, then 
$$
0<z_4\leqslant e^{z_4}-1<\dfrac{2}{\alpha^{n_1}}.
$$

Suppose now that $z_4<0$. We have $2/\alpha^{n_1}<1/2$ for all $n_1>350$. Then, from (\ref{eq:28}), we have 
$$
|1-e^{z_4}|<\dfrac{1}{2}
$$
and therefore $e^{|z_3|}<2$.

Since $z_3<0$, we have : 
$$
0<|z_3|\leqslant e^{|z_3|}-1=e^{|z_3|}\left|e^{|z_3|}-1\right|< \dfrac{4}{\alpha^{n_1}}
$$
which gives 
$$
0<|z_3|<\dfrac{4}{\alpha^{n_1}}
$$
for the both cases ($z_3< 0$ and $z_3>0$ ) and holds for all $n_1>350.$

Replacing the expression of $z_3$ in the above inequality and arguing again as before, we conclude that
\begin{equation}\label{eq:29}
0< \left|a\left(\dfrac{\log 2}{\log \alpha}\right)-n_1+ \dfrac{\log \Upsilon_3(n_1-n_2, n_1-n_3, n_1-n_4)}{\log \alpha}\right|<9\cdot \alpha^{-n_1}.
\end{equation}

Here, we also take, $M:=2.8\times 10^{58}$ and we apply Lemma \ref{lem1} last time in inequality (\ref{eq:29}) for all choices $n_1-n_2\in \{0,314\}$, $n_1-n_3\in \{0, 314\}$ and $n_1-n_4 \in \{0,314\}$ with $(n_1, n_2, n_3, n_4, a)$ a possible solution of equation (\ref{eq:1}), and by omitting study of special cases (because it give a solution presented in Theorem \ref{theo4} ), we get: 

$$n_1< 320.$$
This is false because our assumption that $n_1 > 350$.

This ends the proof of our main theorem. 

\end{proof}
\section{Comments}

In this paper, we found all instances in which a power of two can be expressed as a sum of
four Fibonacci numbers. Given the results obtained, we can make the following conjecture. 

\begin{con}
Consider  the Diophantine equation
\begin{equation}
F_{n_{1}}+F_{n_{2}}+F_{n_{3}}+F_{n_{4}}=p^a, p\geqslant 2, a\geqslant 2 
\end{equation} 
where $n_1, n_2, n_3, n_4, a$  are positive integers with $n_1\geqslant n_2\geqslant n_3\geqslant n_4$ and $p$ is prime,   then $p=2, 3, 5, 7.$ 
\end{con}

\vspace{2cm}
\begin{multicols}{2}
\begin{flushleft}
\begin{center}
\textbf{Pagdame TIEBEKABE}\\
\small{
Université Cheikh Anta Diop (UCAD),}
\footnotesize{Laboratoire d'Algèbre, de Cryptologie,
de Géométrie Algébrique et Applications (LACGAA)}\\ 
\sc{Dakar, Sénégal}
\end{center} 
\end{flushleft}
\columnbreak
\begin{flushright}
\begin{center}
\textbf{Ismaïla DIOUF}\\
\small{
Université Cheikh Anta Diop (UCAD),}
\footnotesize{ Laboratoire d'Algèbre, de Cryptologie,
 de Géométrie Algébrique et Applications (LACGAA)}\\
\sc{Dakar, Sénégal}
\end{center}
\end{flushright}
\end{multicols}

 \end{document}